\documentclass[a4paper,12pt, twoside]{amsart}
\usepackage{amsthm, amsmath}
\usepackage{amssymb,amsfonts,amscd,verbatim, longtable}

\usepackage{ifpdf}
\ifpdf
\usepackage{hyperref}
\else
\usepackage[hypertex]{hyperref}
\fi

\newlength{\myevenmargin}
\DeclareMathOperator\CHAR{char}
\DeclareMathOperator\Bl{Bl}

\def\k{\Bbbk}
\def\bb#1{\mathbb #1}
\def\refeq#1{$(\ref{#1})$}
\def\P{\bb P}

\def\pmat#1{\begin{pmatrix}#1\end{pmatrix}}

\let\le\leq

\let\star *
\let\subset\subseteq

\def\iso{\cong}

\title[On the Poincar\'e polynomials of moduli spaces of one-dimensional sheaves]{An observation on the Poincar\'e polynomials of moduli spaces of one-dimensional sheaves}
\author{Oleksandr Iena}
\email{o.g.yena@gmail.com}
\begin{document}

\begin{abstract}
We notice that for $0<d\le 6$ the Poincar\'e polynomial of Simpson moduli space $M_{dm + 1}(\P_2)$  is divisible by the Poincar\'e polynomial of the projective space $\P_{3d-1}$. 
A somehow regular behaviour of the difference of the  Poincar\'e polynomials of the Hilbert scheme of $\frac{(d-2)(d-1)}{2}$ points on $\P_2$ and the moduli space of Kronecker modules $N(3; d-2, d-1)$ is noticed for $d=4, 5, 6$.  
\end{abstract}
\maketitle
\subsection*{Notations}
 Fix an algebraically closed field $\k$, $\CHAR \k=0$. Let $V$ be a $3$-dimensional vector space over $\k$ and let $\P_2=\P V$ be the corresponding projective plane. Consider a linear polynomial $P(m)=dm +1$  in $m$ with integer coefficients, $d>0$. Let $M_{dm+1}=M_{dm+1}(\P_2)$ be the Simpson moduli  space (cf.~\cite{Simpson1}) of semi-stable sheaves on $\P_2$ with Hilbert polynomial $dm+1$.

\subsection*{Moduli spaces description}
It has been shown in~\cite{LePotier} that 
$M_{dm+1}\iso \P(S^d V^*)$ for $d=1$, $d=2$, and  $M_{3m+ 1}$ is isomorphic to the universal cubic plane curve
\(
\{(C, p)\in \P(S^3V^*)\times \P_2 \mid p\in C\},
\)
which is a $\P_8$-bundle over $\P_2$.

In~\cite{DrezetMaican4m}, \cite{Maican5m}, and~\cite{Maican6m}  the moduli spaces
$M_{dm+1}(\P_2)$ are described in terms of stratifications for $d=4$, $d=5$, and $d=6$  respectively. Similar stratifications are also obtained in~\cite{Yuan}. A description of  $M_{4m+1}(\P_2)$ as a blow-down of a blow-up of a certain projective bundle over a smooth $6$-dimensional base is given in \cite{ChungGlobal} and~\cite{IenaGlobal4}. 
The moduli spaces $M_{dm+1}(\P_2)$ were also studied using wall-crossing techniques  in~\cite{ChoiChung} for $d=4, 5$ and \cite{ChoiChung2} for $d=6$.

\subsection*{Birational models}
As shown in~\cite{MaicanTwoSemiSt}, $M_{dm+1}$ is birational to a $\P_{3d-1}$-bundle over the moduli space  of Kronecker modules $N(3;d-2, d-1)$(cf.~\cite{Drezet, DrezetAltern}). At the same time $M_{dm+1}$ is also birational to the flag Hilbert scheme $H(l, d)$ of pairs $Z\subset C$, where $Z$ is a zero-dimensional scheme of length $l=\frac{(d-2)(d-1)}{2}$ on a planar curve $C\subset \P_2$ of degree $d$.
There is a natural morphism from $H(l, d)$ to the Hilbert scheme $\P_2^{[l]}$ of zero-dimensional subschemes in $\P_2$ of length $l$.
As mentioned in~\cite{ChoiChung}, for $d<6$, $H(l, d)$ is a    
$\P_{3d-1}$-bundle over $\P_2^{[l]}$.

\subsection*{Poincar\'e polynomials}
Clearly, 
\[
P_{M_{m+1}}(t)=P_{\P_2}(t), \quad P_{M_{2m+1}}(t)=P_{\P_5}(t),\quad
P_{M_{3m+1}}(t)=P_{\P_8}(t)\cdot P_{\P_2}(t).
\]
The Poincar\'e polynomials $P_{M_{dm+1}}(t)$ of the moduli spaces $M_{dm+1}(\P_2)$, $d=4,5,6$, have been computed by different authors using different methods. For example, for  $d=4, 5$ the corresponding values can be found in~\cite{ChoiChung}.  $P_{M_{6m+1}}(t)$  is computed in~\cite{ChoiChung2}.
For completeness we provide here the corresponding expressions.
\begin{align*}
P_{M_{4m+1}}(t)=&1+2t^2+6t^4+10t^6+14t^8+15t^{10}+16t^{12}+16t^{14}+16t^{16}+\\
&16t^{18}+16t^{20}+16t^{22}+15t^{24}+14t^{26}+10t^{28}+6t^{30}+2t^{32}+t^{34},\\
\end{align*}
\begin{align*}
P_{M_{5m+1}}(t)=&1+2t^2+6t^4+13t^6+26t^8+45t^{10}+68t^{12}+87t^{14}+100t^{16}+\\
&107t^{18}+111t^{20}+112t^{22}+113t^{24}+113t^{26}+113t^{28}+112t^{30}+111t^{32}+107t^{34}+\\
&100t^{36}+87t^{38}+68t^{40}+45t^{42}+26t^{44}+13t^{46}+6t^{48}+2t^{50}+t^{52},\\
\end{align*}
\begin{align*}
P_{M_{6m+1}}(t)=&(1+t^2+4t^4+7t^6+16t^8+25t^{10}+47t^{12}+68t^{14}+104t^{16}+128t^{18}+146t^{20}+\\
&128t^{22}+104t^{24}+
68t^{26}+47t^{28}+25t^{30}+16t^{32}+7t^{34}+4t^{36}+t^{38}+t^{40}
)\cdot\frac{t^{36}-1}{t^2-1}.
\end{align*}

\subsection*{Observing a regular behaviour}
Notice that the expression for $P_{M_{6m+1}}(t)$ in~\cite{ChoiChung2} is given as a multiple of polynomial
\[
\frac{t^{36}-1}{t^2-1}=1+t^2+\dots+t^{32}+t^{34},
\]
which is the Poincar\'e polynomial of the projective space $\P_{17}$. 

Decomposing the polynomials  $P_{M_{4m+1}}(t)$ and  $P_{M_{5m+1}}(t)$ into irreducible factors using~\textsc{Singular}~\cite{SingularProgram}, we notice that for every $0<d\le 6$ the Poincar\'e polynomial $P_{M_{dm+1}}(t)$ is divisible by the Poincar\'e polynomial of the projective space $\P_{3d-1}(t)$, i.~e.,   
 $P_{M_{dm+1}}(t)$ looks for $0<d\le 6$ as the Poincar\'e polynomial of a projective $\P_{3d-1}$-bundle over some space.
 
Denote 
 \[
 P_{v, d}(t)=\cfrac{P_{M_{dm+1}}(t)}{P_{\P_{3d-1}}(t)},\quad 0<d \le 6.
 \] 
 Then
 \begin{align*}
 P_{v, 1}(t)=&1,\quad
 P_{v, 2}(t)=1,\quad
 P_{v, 3}(t)=1+t^2+t^4,\\
 P_{v, 4}(t)=&t^{12}+t^{10}+4t^8+4t^6+4t^4+t^2+1,\\
 P_{v, 5}(t)=&t^{24}+t^{22}+4t^{20}+7t^{18}+13t^{16}+19t^{14}+23t^{12}+19t^{10}+13t^8+7t^6+4t^4+t^2+1\\
 P_{v, 6}(t)=&t^{40}+t^{38}+4t^{36}+7t^{34}+16t^{32}+25t^{30}+47t^{28}+68t^{26}+104t^{24}+128t^{22}+146t^{20}+\\
 &128t^{18}+104t^{16}+68t^{14}+47t^{12}+25t^{10}+16t^8+7t^6+4t^4+t^2+1.
 \end{align*}

 Computing the Poincar\'e polynomials of $N(3; d-2, d-1)$ and $\P_2^{[l]}$ using the formulas from~\cite{DrezetCohom} and~\cite{GoettscheFormula} and their computer algebra implementations in~\cite{GoettscheLib}, one notices that all the coefficients of $P_{v, d}(t)$, $d=3, 4, 5, 6$, are between the values of the corresponding Betti numbers of the moduli space of Kronecker modules $N(3; d-2, d-1)$ and the Hilbert scheme $\P_2^{[l]}$. More precisely,
 \begin{align*}
 P_{v, 3}(t)=&1+t^2+t^4=P_{N(3;1,2)}(t)=P_{\P_2^{[1]}}(t),\\
 &\\
 P_{v, 4}(t)=&P_{N(3;2,3)}(t)+t^4(1+t^2+t^4)=P_{\P_2^{[3]}}(t)-t^2(t^4+1)(t^4+t^2+1),\\
 &\\
 P_{v, 5}(t)=&P_{N(3;3,4)}(t)+t^4(1+t^2+t^4)(1+t^2+3t^4+5t^6+3t^8+t^{10}+t^{12})=\\
             & P_{\P_2^{[6]}}(t)- t^2(t^2+1)^2(t^{16}+5t^{12}+3t^{10}+9t^8+3t^6+5t^4+1), \\
 &\\
 P_{v, 6}(t)=&P_{N(3;4,5)}(t)+t^4(1+t^2+t^4)\cdot f=P_{\P_2^{[10]}}(t)-t^2(t^4+t^2+1)\cdot g,
 \end{align*}
 where
 \begin{align*}
  f=&1+t^2+4t^{4}+6t^{6}+14t^{8}+18t^{10}+31t^{12}+33t^{14}+\\
  &31t^{16}+18t^{18}+14t^{20}+6t^{22}+4t^{24}+t^{26}+t^{28},\\
 g=&1+t^{2}+4t^{4}+8t^{6}+20t^{8}+35t^{10}+66t^{12}+93t^{14}+108t^{16}+\\
 &93t^{18}+66t^{20}+35t^{22}+20t^{24}+8t^{26}+4t^{28}+t^{30}+t^{32}.
 \end{align*}

\subsection*{Questions to answer}
We formulate  here some questions that seem reasonable to ask.
\begin{enumerate}
\item Is it a coincidence that $P_{M_{dm+1}}(t)$ is divisible by $P_{\P_{3d-1}}(t)$ for $0<d \le 6$?

\item Can one expect this also to be the case for $d>6$?

\item Are there meaningful geometric spaces with Poincar\'e polynomials $P_{v, d}(t)$? 

\end{enumerate}

\subsection*{Remarks on the Poincar\'e polynomials of Hilberts schemes of points and moduli spaces of Kronecker modules}
As a somehow related side remark we share here some observations on 
the difference of the  Poincar\'e polynomials of the Hilbert scheme of $l$ points on $\P_2$ and the moduli space of Kronecker modules $N(3; d-2, d-1)$.

Notice that the schemes $\P_2^{[l]}$ and $N(3; d-2, d-1)$ are birational. 
This can be explained as follows. Let $H_d'\subset H_d$ be the closed subscheme of schemes lying on a curve of degree $d-3$. It is an irreducible hypersurface in $H_d$. Let $N'_d\subset N_d$ be the closed subscheme consisting of the classes of Kronecker modules whose maximal minors have a common factor. Then $N_d\setminus N'_d$ is isomorphic to $H_d\setminus N'_d$, the isomorphism sends the a class of a Kronecker module to the vanishing scheme of its maximal minors.

For  $d=3$ there is clearly an isomorphism $\P_2^{[1]}\iso N(3; 1, 2)$. 
For $d=4$ the Hilbert scheme $\P_2^{[3]}$ is a blow-up  of $N(3; 2, 3)$ along a smooth subscheme that is isomorphic to a projective plane 
(cf.~\cite[Th\'eor\`eme~4]{DrezetCohom}). Though the explicit description of this birational equivalence is unknown to the author for $d>4$, we wish to provide here the following observations.

First of all consider the difference
\begin{align}\label{eq: difference d=4}
\begin{split}
P_{\P_2^{[3]}}(t)-P_{N(3; 2, 3)}(t)=&t^2(1+t^2+t^4)^2=\\
&(1+t^2+t^4)(1+t^2+t^4+t^6-1)=\\
&P_{\P_2}(t)(P_{\P_3}(t)-1)=P_{\P_2}(t)P_{\P_3}(t)-P_{\P_2}(t),
\end{split}
\end{align}
which indeed reflects the fact that $\P_2^{[3]}$ is obtained from $N(3; 2, 3)$ by a substitution of a subvariety isomorphic to a projective plane by a $\P_3$-bundle over it.

At the same time the differences
\begin{align} \label{eq: difference for d=5}
P_{\P_2^{[6]}}(t)-P_{N(3; 3, 4)}(t)=t^2(1+t^2+t^4)^2(1+t^2+3t^4+7t^6+3t^8+t^{10}+t^{12})
\end{align}
and 
\begin{align} \label{eq: difference for d=6}
&P_{\P_2^{[10]}}(t)-P_{N(3; 4, 5)}(t)= t^2(1+t^2+t^4)^2\cdot f,
\end{align}
with
\(f=1+t^2+3t^{4}+8t^{6}+15t^{8}+26t^{10}+43t^{12}+55t^{14}+43t^{16}+26t^{18}+15t^{20}+8t^{22}+3t^{24}+t^{26}+t^{28}\), surprisingly turn out to be multiples of~\refeq{eq: difference d=4}.

Concerning~\refeq{eq: difference for d=5}, one can easily notice that $N'_5$ contains a closed subvariety $N''$ that corresponds to the Kronecker modules with maximal minors having a common quadratic factor $q$. The corresponding points are the equivalence classes of the Kronecker modules
\[
\pmat{
0&x_2&-x_1&l_0\\
-x_2&0&x_0&l_1\\
x_1&-x_0&0&l_2
}
\]
such that $q=l_0x_0+l_1x_1+l_2x_2$. Here $x_0, x_1, x_2$ is a fixed basis of $V^*$. Then $N''$ is isomorphic to the space of conics, i.~e., $N''\iso \P_5$. Then
\begin{align*}
P_{\Bl_{N''}N_5}(t)-P_{N_5}(t)=
&P_{\P_5}(t)P_{\P_6}(t)-P_{\P_5}(t)=P_{\P_5}(P_{\P_6}-1)=t^2P_{\P_5}^2=\\
&t^2P_{\P_2}^2(1+t^6)^2=
t^2P_{\P_2}^2(1+2t^6+t^{12})=t^2(1+t^2+t^4)^2(1+2t^6+t^{12})
\end{align*}
because
\[
P_{\P_5}=\frac{1-t^{12}}{1-t^2}=\frac{1-t^{6}}{1-t^2}(1+t^6)=P_{\P_2}\cdot(1+t^6).
\]
So, indeed, $H_5$ seems to be not so far away from being the blow-up of $N_5$ along $N''$.

One could also expect~\refeq{eq: difference for d=6} to bear some resemblances with the difference
\begin{align*}
P_{\P_9}(t)P_{\P_{10}}(t)-P_{\P_9}(t)=P_{\P_9}(t)(P_{\P_{10}}(t)-1)=t^2P_{\P_9}(t)^2
\end{align*}
corresponding to a blow-up of $N_6$ at a subvariety isomorphic to the space of cubic planar curves. In this case, however, the factor $(1+t^2+t^4)^2$ does not appear immediately as a factor of $t^2P_{\P_9}(t)^2$.

One easily checks using~\cite{GoettscheLib} that for $d>6$ the differences $P_{\P_2^{[l]}}(t)-P_{N(3; d-2, d-1)}(t)$ are not divisible by~\refeq{eq: difference d=4}.


\begin{thebibliography}{10}

\bibitem{ChoiChung2}
Jinwon {Choi} and Kiryong {Chung}.
\newblock {The geometry of the moduli space of one-dimensional sheaves.}
\newblock {\em {Sci. China, Math.}}, 58(3):487--500, 2015.

\bibitem{ChoiChung}
Jinwon {Choi} and Kiryong {Chung}.
\newblock {Moduli spaces of $\alpha$-stable pairs and wall-crossing on
  $\mathbb{P}^2$.}
\newblock {\em {J. Math. Soc. Japan}}, 68(2):685--709, 2016.

\bibitem{ChungGlobal}
K.~{Chung} and H.-B. {Moon}.
\newblock {Chow ring of the moduli space of stable sheaves supported on quartic
  curves}.
\newblock {\em Q. J. Math}, 68(3):851--887, September 2017.

\bibitem{SingularProgram}
Wolfram Decker, Gert-Martin Greuel, Gerhard Pfister, and Hans Sch\"onemann.
\newblock {\sc Singular} {4-0-2} --- {A} computer algebra system for polynomial
  computations.
\newblock \url{http://www.singular.uni-kl.de}, 2015.

\bibitem{Drezet}
J.-M. Drezet.
\newblock {Fibr\'es exceptionnels et vari\'et\'es de modules de faisceaux semi-
  stables sur $\mathbb{P}\sb 2(\mathbb{C})$. (Exceptional bundles and moduli
  varieties of semi-stable sheaves on $\mathbb{P}\sb 2(\mathbb{C}))$.}
\newblock {\em J. Reine Angew. Math.}, 380:14--58, 1987.

\bibitem{DrezetCohom}
J.-M. {Drezet}.
\newblock {Cohomologie des vari\'et\'es de modules de hauteur nulle.
  (Cohomology of moduli varieties of height zero).}
\newblock {\em {Math. Ann.}}, 281(1):43--85, 1988.

\bibitem{DrezetAltern}
Jean-Marc Dr{\'e}zet.
\newblock Vari\'et\'es de modules alternatives.
\newblock {\em Ann. Inst. Fourier (Grenoble)}, 49(1):v--vi, ix, 57--139, 1999.

\bibitem{DrezetMaican4m}
Jean-Marc Dr{\'e}zet and Mario Maican.
\newblock On the geometry of the moduli spaces of semi-stable sheaves supported
  on plane quartics.
\newblock {\em Geom. Dedicata}, 152:17--49, 2011.

\bibitem{GoettscheFormula}
Lothar {G\"ottsche}.
\newblock {Betti numbers for the Hilbert function strata of the punctual
  Hilbert scheme in two variables.}
\newblock {\em {Manuscr. Math.}}, 66(3):253--259, 1990.

\bibitem{GoettscheLib}
Oleksandr Iena.
\newblock {\sc goettsche.lib}, a {\sc {s}ingular} library implementing some
  formulas for {B}etti numbers (by {D}rez\'et, {G}\"ottsche, {N}akajima and
  {Y}oshioka, {M}acdonald).
\newblock
  \url{https://github.com/Singular/Sources/blob/spielwiese/Singular/LIB/goettsche.lib},
  2016--2018.

\bibitem{IenaGlobal4}
Oleksandr {Iena}.
\newblock {On the fine Simpson moduli spaces of $1$-dimensional sheaves
  supported on plane quartics}.
\newblock {\em Open Math.}, 16(1):46--62, 2018.

\bibitem{LePotier}
J.~Le~Potier.
\newblock Faisceaux semi-stables de dimension {$1$} sur le plan projectif.
\newblock {\em Rev. Roumaine Math. Pures Appl.}, 38(7-8):635--678, 1993.

\bibitem{MaicanTwoSemiSt}
Mario Maican.
\newblock {On two notions of semistability.}
\newblock {\em Pac. J. Math.}, 234(1):69--135, 2008.

\bibitem{Maican5m}
Mario Maican.
\newblock On the moduli spaces of semi-stable plane sheaves of dimension one
  and multiplicity five.
\newblock {\em Illinois J. Math.}, 55(4):1467--1532 (2013), 2011.

\bibitem{Maican6m}
Mario Maican.
\newblock The classification of semistable plane sheaves supported on sextic
  curves.
\newblock {\em Kyoto J. Math.}, 53(4):739--786, 2013.

\bibitem{Simpson1}
Carlos~T. Simpson.
\newblock Moduli of representations of the fundamental group of a smooth
  projective variety. {I}.
\newblock {\em Inst. Hautes \'Etudes Sci. Publ. Math.}, (79):47--129, 1994.

\bibitem{Yuan}
Yao Yuan.
\newblock Moduli spaces of semistable sheaves of dimension 1 on {$\mathbb
  P^2$}.
\newblock {\em Pure Appl. Math. Q.}, 10(4):723--766, 2014.

\end{thebibliography}

\def\cprime{$'$} \def\cprime{$'$} \def\cprime{$'$}

\end{document}